\documentclass[12pt]{article}

\RequirePackage[OT1]{fontenc}
\RequirePackage{amsthm,amsmath}
\RequirePackage{amssymb}
\RequirePackage[colorlinks,citecolor=blue,urlcolor=blue]{hyperref}
\usepackage{tabularx}\newcolumntype{C}{>{\centering\arraybackslash}X}

\usepackage{longtable}
\usepackage[utf8]{inputenc}
\usepackage{lineno}
\usepackage{float}

\usepackage{hyperref}

\usepackage{stackrel}

\usepackage{makecell}

\usepackage[round]{natbib}
\usepackage{multirow}
\usepackage[active]{srcltx}
\usepackage{graphicx}
\usepackage{verbatim,enumerate}
\usepackage{rotating, lscape}
\usepackage[usenames,dvipsnames]{color}
\usepackage{subcaption}

\usepackage{enumerate}

\usepackage{soul}
\usepackage[makeroom]{cancel}

\def\S{\mathbb{S}}

\newcommand{\R}{\mathbb{R}}
\newcommand{\Z}{\mathbb{Z}_+}
\numberwithin{equation}{section}
\theoremstyle{plain}
\newtheorem{thm}{Theorem}[section]
\newtheorem{cor}[thm]{Corollary}

\theoremstyle{definition}

\newtheorem{rem}[thm]{Remark}

\begin{document}

	\fontsize{12}{14pt plus.8pt minus .6pt}\selectfont \vspace{0.8pc}
	\centerline{\LARGE \bf Dimension Walks on Generalized Spaces \\ \\
	}
	
	\vspace{.8cm} 
	\begin{center}

		{\large {\sc Marcos Lopez de Prado}}
		\footnote{\baselineskip=10pt
			ADIA $\&$ Cornell $\&$ KU \\
			E-mail: mldp@truepositive.com 
		}

		{\large {\sc Ana Paula Peron}}
		\footnote{\baselineskip=10pt
			Department of Mathematics,\\
			University of S\~ao Paulo, S\~ao Carlos, Brazil.
			\\ 
			E-mail: apperon@icmc.usp.br  \\
			Partially supported by Fapesp, Brazil.
		}

		{\large {\sc Emilio Porcu}}
		\footnote{\baselineskip=10pt
				Department of Mathematics, \\
				Khalifa University, Abu Dhabi, United Arab Emirates, \\
			$\&$ School of Computer Science and Statistics, Trinity College Dublin.
			\\ 
			E-mail: emilio.porcu@ku.ac.ae  
		}\footnote{Corresponding author}
		
	\end{center}

	\centerline{\it This paper is dedicated to Daryl J. Daley and Robert Schaback. Such beautiful minds.}
	
	{
		\hypersetup{linkcolor=black}
		%		\tableofcontents
	}

%    Abstract is required.
\begin{abstract}
Let $d,k$ be positive integers. We call generalized spaces the cartesian product of the $d$-dimensional sphere, $\S^d$, with the $k$-dimensional Euclidean space, $\R^k$.
We consider the class ${\mathcal P}(\S^d \times \R^k)$ of continuous functions $\varphi: [-1,1] \times [0,\infty) \to \R$ such that the mapping $ C: \left ( \S^d \times \R^k \right )^2 \to \R$, defined as $C \Big ( (x,y),(x^{\prime},y^{\prime})\Big ) = \varphi \Big ( \cos \theta(x,x^{\prime}), \|y-y^{\prime}\| \Big )$, $(x,y), \; (x^{\prime},y^{\prime}) \in \S^d \times \R^k$, is positive definite. We propose linear operators that allow for walks through dimension within generalized spaces while preserving positive definiteness. 
\end{abstract}


%    Text of article.
\section{Introduction}
\subsection{Context} 
The paper deals with the problem of linear projection operators that map the set of positive definite functions on a given space into the set of positive definite functions on lower or higher dimensional space. Specifically, we consider continuous functions that are positive definite over generalized spaces, that we define here as the Cartesian product of the $k$-dimensional Euclidean space with a $d$-dimensional unit sphere. \\
There are several motivations to consider this setting, and we sketch them as follows: 

a) There has been an increasing interest from several branches of statistics, machine learning, and finance, for positive definite functions defined over these product spaces, and the reader is referred to the recent review by \cite{porcu30}. The applications to real cases are ubiquitous, ranging from climate and atmospheric sciences until deep learning on manifolds. \\
As far as finance is concerned, Gaussian processes play a central role in financial modelling. The field of econometrics has devoted much effort to the modelling of financial time series \citep{hamilton}. The Brownian motion is essential to the pricing of financial derivatives \citep{hull}. The Ornstein-Uhlenbeck process is often used to develop investment and trading strategies \citep{LipLop}. Gaussian processes are one-dimensional applications of the general concept of Gaussian random field (GRF). We motivate some of the uses of GRFs in finance.
A key feature of financial datasets is time and spatial dependence. Coetaneous observations from variables in close proximity tend to be more similar. For example, returns from U.S. stocks last month are more similar to returns from U.S. stocks this month than returns from U.S. stocks one year ago, or returns from Chinese stocks last month. \cite{willinger} noted that a Brownian motion with drift does not replicate the time dependence observed in asset returns. Not surprisingly, GRFs have attracted considerable interest among researchers interested in modelling the joint time-space dynamics of financial processes. 
To cite a few examples, \cite{kennedy} and \cite{goldstein} modelled the term structure of interest rates as a two-dimensional random field. In their models, time increments are independent, while the correlation structure between bond yields of different maturities can be modelled with great flexibility. \cite{kimmel} enhanced this approach by adding a state-dependent volatility. \cite{albeverio} introduced L{\'e}vy fields to the modelling of yield curves. \cite{ozkan} applied random fields to incorporate credit risk to the modelling of yield curves.
As important as the term structure of interest rates is, it is not the only financial application of Gaussian fields. At least two further applications stand out: Option pricing and actuarial modelling. For example, \cite{hainaut} proposes an alternative model for asset prices with sub-exponential, exponential and hyper-exponential autocovariance structures. Hainaut sees price processes as conditional Gaussian fields indexed by the time. Under this framework, option prices can be computed using the technique of the change of numeraire. \cite{biffis} applied random fields to modelling the intensity of mortality, in an attempt to incorporate cross-generation effects. \cite{biagini} built on that work to price and hedge life insurance liabilities.

b) Several branches of spatial statistics and computer sciences are interested in the simulation of random processes defined over generalized spaces, and we refer the reader to \cite{emery2008turning}. It turns out that the use of these operators becomes crucial when associated with turning bands techniques \citep{matheron1973intrinsic}, which allow for simulation of a given random process from projections on lower dimensional spaces. 

c) Projection operators for radial positive definite functions allowed to build positive definite functions that are compactly supported over balls embedded in $k$-dimensional Euclidean spaces. This inspired a fertile literature from spatial statistics with the goal of achieving accurate estimates while allowing for computational scalability. For instance, the tapering approach \citep{furrer2006covariance} is substantially based on this idea.  

d) There is a fertile literature from projection operators for symmetric (or radially symmetric) distributions, where radial symmetry is intended with respect to the composition of a given candidate function with the classical $\alpha$-norms \citep{cambanis1983alpha}.

\subsection{Literature Review}

Let $\R^k$ denote the $k$-dimensional Euclidean space, and let $\S^d$ be the $d$-dimensional unit sphere embedded in $\R^{d+1}$. Let $\|\cdot\|$ denote Euclidean distance and $\theta(x,y): = \arccos \left ( \langle x,y \rangle \right ) $ denote the geodesic distance in $\S^d$, with $\langle \cdot, \cdot \rangle$ denoting the dot product in $\R^{d+1}$. A continuous function $C: \R^k \to \R$ is called radially symmetric if there exists a continuous function $f: [0,\infty) \to \R$ such that $C(x)=f \circ \|x\|$, $x \in \R^k$, with $\circ$ denoting composition. The function $f$ is called the radial part of $C$. Radial symmetry is known as {\em isotropy} in spatial statistics \citep{daley-porcu}. 
A function $C: \S^d \times \S^d \to \R$ is called geodesically isotropic if $C(x,y)= g \circ \theta(x,y) $ for some continuos function $g:[0,\pi] \to \R$. 

 Positive definite functions that are radially symmetric over $k$-dimensional Euclidean spaces have a long history that can be traced back to \cite{scho38}. Projections operators that map a positive definite radial mapping from $\R^k$ into $\R^{k \pm h}$, for $h$ a positive integer, have been considered in the Matheron's clavier spherique \citep{matheron1973intrinsic, matheron1970random}.  Matheron coined the terms {\em descente} and {\em mont{\'e}e} to define special operators that will be described throughout. The terms originate from an appealing physical interpretation in a mining context.  These projections operators have then be investigated by \cite{eaton1981projections}, and subsequently by \cite{wendland1995piecewise}, \cite{schaback1996operators}, and by \cite{gneiting2002} in the context of positive definite radial functions that are additionally compactly supported on balls embedded in $\R^k$ with given radii. The work by \cite{daley-porcu} provides a general perspective of such operators, in concert with some generalizations of the previously mentioned works. These linear operators have turned to be very useful to establish criteria of the P{\'o}lya type for radially symmetric positive definite functions \citep{gneiting2001criteria}, as well as in the definition of multiradial positive definite functions \citep{porcu2007descente}. In probability theory, similar projection operators turned useful in the seminal paper by \cite{cambanis1983alpha} and in \cite{fang2018symmetric}. 
 
 Positive definite functions that are geodesically isotropic on $d$-dimensional spheres have been characterized in \cite{schoenberg1942}. Projection operators for this class of functions have been studied to a limited extent only, and we refer to the recent papers by \cite{beatson2014polya} and more recently to the same authors \citep{beat-zucatell-2016, beatson2017dimension}. Properties of these operators have then been inspected in \cite{trubner2017derivatives}.

\subsection{The Problem, and Our Contribution}

A characterization of projection operators on product spaces of the type $\S^d \times \R^k$ has been elusive so far. The only exception being \cite{bing-sym2019}, who consider the product space
 $\S^d \times \R$, and projections that are defined marginally for the sphere only. 

Our paper contributes to the literature as follows. 
In Section 2 we provide the notations and basic literature. In Section 3 we define the Descente and Montée operators on the  generalized space $\S^d\times\R^k$. The main results are statement in Section 4 and their proofs in Appendix A.

\section{Notations and Background}

Let $X,Y$ be nonempty sets. A function $C: (X \times Y)^2 \to \R$ is called positive definite if, for any finite system $\{ a_k \}_{k=1}^N\subset\R$ and points $\{ (x_k,y_k) \}_{k=1}^N\subset X\times Y$, the following inequality is preserved: 
$$ \sum_{k=1}^N \sum_{h=1}^N a_k C \left ( (x_k,y_k),(x_h,y_h) \right ) a_h \ge 0. $$
We deal with the case $X= \S^d$ and $Y= \R^k$, for $d$ and $k$ being positive integers. Additionally, we suppose $C$ to be continuous, and that there exists a continuous function $\varphi: [-1,1] \times [0,\infty) \to \R$ such that 
\begin{equation}
\label{eugenio1} C \Big ( (x,y),(x^{\prime},y^{\prime})\Big ) = \varphi \Big ( \cos \theta(x,x^{\prime}), \|y-y^{\prime}\| \Big ), \qquad (x,y), \; (x^{\prime},y^{\prime}) \in \S^d \times \R^k.  
\end{equation}
We call ${\mathcal P}(\S^d \times \R^k)$ the class of such functions, $\varphi$. Analogously, we call ${\mathcal P}(\S^d)$ the class of continuous functions $\psi: [-1,1] \to \R$ such that for the function $C  $ in (\ref{eugenio1}) it is true that, for $y=y^{\prime}$, $C \left ( (x,y),(x^{\prime},y) \right ) = \varphi(\cos \theta(x,x^{\prime}),0) = \psi(\cos \theta(x,x^{\prime}))$. The class ${\mathcal P}(\R^k)$ is defined analogously. The classes ${\mathcal P}(\R^k)$ and ${\mathcal P}(\S^d)$  have been characterized by  \cite{scho38} and \cite{schoenberg1942}, respectively. 
The class ${\mathcal P}(\S^d \times \R^k)$ has been characterized by \cite{porcu-berg} through a uniquely determined expansions of the type 
	\begin{equation}\label{eq-pd_SdxRk}
	\begin{array}{cc}
\displaystyle	\varphi(x,t)=\sum_{n=0}^{\infty}f_{n}^d(t)C_{n}^{(d-1)/2}(x), \quad  (x,t)\in[-1,1]\times[0,\infty),\\
\text{where the functions $f_n^d$ belong to $\mathcal P(\R^k)$, $n\in\Z$, and $\displaystyle\sum_{n=0}^{\infty}f_{n}^d(0){C_{n}^{(d-1)/2}(1)}<\infty$.}
	\end{array}
\end{equation}	
The expansion above is uniformly convergent on $[-1,1]\times[0,\infty)$. The coefficients functions $f_n^d$ are called $d$-Schoenberg functions of $\varphi$. The functions $C_{n}^{(d-1)/2}$ are the %\rosa{\st{normalized}} 
Gegenbauer polynomials of degree $n$ associated to the index $(d-1)/2$ \citep{szego}.%\rosa{\st{, being identically equal to $1$ for $x=1$}}.

Proposition 3.8 in \cite{porcu-berg} shows that if $\varphi$ belongs to the class $\mathcal P(\S^{d}\times\R^k)$, then is continuously differentiable with respect to the first variable.

{
	It is also important to note that a continuous function $x\in[-1,1]\mapsto \varphi(x,t)$ has an Abel-summable expansion  for each $t\in[0,\infty)$  in the form \citep[see the proof of Theorem 3.3 in][]{porcu-berg}
	\begin{equation}\label{eq-los_grandes_BP}
		\varphi(x,t) \sim \sum_{n=0}^{\infty}f_{n}^d(t)C_{n}^{(d-1)/2}(x),
	\end{equation} 
	where 
	\begin{equation}\label{eq-coef-los_grandes_BP}
		f_{n}^d(t) = \varsigma_n^d \int_{-1}^1\varphi(x,t)C_n^{(d-1)/2}(x)(1-x^2)^{d/2-1}dx,
	\end{equation}
	and $\varsigma_n^d$ are positive constants.
}

\subsection{Some Useful Facts}

Arguments in \cite{scho38} prove that, for every $n=0,1,\ldots$, each function $f_n^d\in\mathcal P(\R^k)$ in (\ref{eq-pd_SdxRk}) admits a  uniquely determined Riemann-Stieltjes integral representation of the form
\begin{equation}\label{eq-pd-Rk}
f_{n}^d(t) =\int_0^\infty\Omega_k(tr)dF_n(r), \quad t\in[0,\infty),
\end{equation}
where $F_n$ is a non negative bounded measure on $[0,\infty)$. The function $\Omega_k:[0,\infty)\to\R$ is given by
\begin{equation}\label{eq-Omega}
\Omega_k(t) = \Gamma\left(\frac{k}2\right) \left(\frac2t\right)^{(k-2)/2}J_{(k-2)/2}(t),
\end{equation}
where $J_\nu$ is the Bessel function of the first kind of order $\nu$ given by
$$
J_\nu(t)=\left({\frac {t}{2}}\right)^\nu\sum_{m=0}^{\infty}\frac{(-1)^{m}}{m!\Gamma(m+\nu +1)} \left({\frac {t}{2}}\right)^{2m}.
$$
We follow \cite{daley-porcu} and we call  $F_n$ the $k$-Schoenberg  measure of $f_n^d$. We also note that we are abusing of notation when writing $F_{n}$ instead of $F_{n}^d$. This last notation will not be used unless explicitly needed. 

Some technicalities will be exposed here to allow for a neater exposition. 
The derivative function of the function $\Omega_k$ is uniformly bounded, and it is given by \citep[see ][]{daley-porcu, gneiting2002,porcu2007descente}
\begin{equation}\label{eq-derv-Omega}
\frac{d\Omega_k}{dt}(t) = \Omega_k'(t) = -\frac1k t \Omega_{k+2}(t), \qquad t \ge 0.
\end{equation}
Also ,
\begin{equation}\label{eq-Omega_ltda}
|\Omega_k(t)|<1 =\Omega_k(0), \quad t>0.
\end{equation}	
Since $\lim_{t\to\infty}\Omega_k(t)=0$ for $k>0$  (see \cite{daley-porcu}), we have
\begin{equation}\label{eq-int-Omega}
\int_t^\infty u \Omega_k(u)du = (k-2)\Omega_{k-2}(t), \quad t\geq0.
\end{equation}
Some properties of Gegenbauer polynomials will turn to be useful throughout. 
For instance, we can invoke $4.7.14$ in  \cite{szego} to infer that  
\begin{equation}\label{eq-derv-Gegenb}
\frac{dC_{n}^{\lambda}}{dx}(x) = (C_{n}^{\lambda})'(x) = 
\delta_\lambda C_{n-1}^{\lambda+1}(x), \qquad -1 \le x \le 1,
\end{equation}
%\ana{I think we have a problem: the formula (2.8) holds for the Gegenbauer polynomials NOT normalized.... Then in (2.2) we can not use the normalized polynomials and them we need $\sum f_n^dC_n^{(d-1)/2}(1)<\infty$..... Unfortunately  I think we need to check the (uniform) convergence of the series in all our results in order to conclude that in each case the operators are positive definite.... I am sorry... I am working on it...}
and, as a consequence
\begin{equation}\label{eq-integ-Gegenb}
\int_{-1}^xC_{n}^{\lambda}(x)dx = 
\frac1{\delta_\lambda} (C_{n+1}^{\lambda-1}(x)-C_{n+1}^{\lambda-1}(-1)), 
\end{equation}
where 
\begin{equation}\label{eq-delta_Lamb}
\delta_\lambda = \left\{\begin{array}{ll} 
	2\lambda, &\lambda> -1/2\,\, (\lambda\neq0),\\\\
	2, &\lambda= 0.
\end{array}
\right.
\end{equation}
Theorem 7.32.1 and Equation 4.7.3 in \cite{szego} show that, for $\lambda>-1/2$, 
\begin{equation}\label{eq-Gegenb-ltda}
|C_{n}^{\lambda}(x)| \leq C_{n}^{\lambda}(1) %=\binom{n+2\lambda-1}{n} 
 =\frac{\Gamma(n+2\lambda)}{\Gamma(n+1)\Gamma(2\lambda)}, \quad x\in[-1,1].
\end{equation}
%\end{eqnarray}

Also, it is true that 
\begin{eqnarray}\label{eq-razao-gegenb}
\frac{C_{n+j}^{\lambda-k}(1)}{C_{n}^{\lambda}(1)} 
&\leq& \frac{\Gamma(2\lambda)}{\Gamma(2\lambda-2k)}:=\varrho_{\lambda,k}, \quad \forall n\in\Z. 
%\nonumber
\end{eqnarray}

The following inequality \citep[see][]{Berg} will be repeatedly used in the manuscript:
$$
|f(t)|\leq f(0),\quad t\in[0,\infty),\quad f\in\mathcal P(\R^k).
$$

We will also make use of the following fact: if  $\varphi:[-1,1]\times\R\to\R$ has a derivative $\varphi_x$  with respect to the first variable for each $t\in[0,\infty)$ and if both functions have Gegenbauer expansions of the form
	\begin{equation}\label{eq-phi__phi_x}
		\varphi_x(x,t) \sim \sum_{n=0}^{\infty} f_{n}^\lambda(t)C_{n}^{\lambda}(x),
		\qquad
		\varphi_x(x,t) \sim \sum_{n=0}^{\infty} \tilde f_{n}^{\lambda+1}(t)C_{n}^{\lambda+1}(x),
	\end{equation} 
	$(x,t ) \in [-1,1 ] \times [0,\infty)$, then
	\begin{equation}\label{eq-coef-phi__phi_x}
		\tilde f_{n-1}^{\lambda+1}(t) = \delta_\lambda f_{n}^\lambda(t), \quad n\in\Z^*, \quad \lambda>0.
	\end{equation}
	The proof is very similar to the proof of Lemma 2.4 in \cite{beat-zucatell-2016} and we omit it for the sake of brevity.

\section{An Historical Account on {\em Mont{\'e}e} and {\em Descente} Operators}

\cite{beatson2017dimension} defined the Descente and Mont{\'e}e operators for the class ${\mathcal P}(\S^{d})$. Specifically, the Descente ${\mathcal D}$ is defined as 
$$  \left ( {\mathcal D} f \right ) (x) = \frac{d}{d x } f(x) = f^{\prime } (x), \qquad x \in [-1,1]
, $$ provided such a derivative exists. The Mont{\'e}e ${\mathcal I}$ is instead defined as 
$$  \left ( {\mathcal I} f \right ) (x) = \int_{-1}^x f(u) d u, \qquad x \in [-1,1]. $$
 \cite{beatson2017dimension} have shown that $f \in {\mathcal P}(\S^{d+2})$ implies that there exists a constant, $\kappa$, such that  $\kappa+ {\mathcal I} f \in {\mathcal P}(\S^{d})$. Also, $f \in {\mathcal P}(\S^{d})$ implies ${\mathcal D} f \in {\mathcal P}(\S^{d+2})$. The implication in terms of differentiability at $x=1$ are nicely summarized therein. \\
 The {\em tour de force} by \cite{beatson2017dimension} has then been generalized by \cite{bing-sym2019}: let $d\in\mathbb N$ and $\varphi:[-1,1]\times\R\to\R$ be a continuous  functions.
The Montée $\mathcal{I}$ and Descente $\mathcal{D}$ operators are defined respectively by
\begin{equation}\label{eq-montee-SXR}
\mathcal{I}(\varphi)(x,t) :=\int_{-1}^x\varphi(u,t)du,,\qquad (x,t)\in[-1,1]\times [0,\infty)
\end{equation}
when $f$ is integrable with respect to the first variable, and
\begin{equation}\label{eq-descente-SXR}
\mathcal{D}(\varphi)(x,t) :=\frac{\partial \varphi}{\partial x}(x,t),\qquad (x,t)\in[-1,1]\times [0,\infty).
\end{equation}
They prove that if  $\varphi\in\mathcal P(\S^{d}\times\R)$, then $\mathcal D\varphi\in\mathcal P(\S^{d+2}\times \R)$
and in their correction of Theorem 2.1 they provided conditions under $\varphi\in\mathcal P(\S^{d+2}\times\R)$ such that $\mathcal I \varphi\in\mathcal P(\S^{d}\times\R)$.

Mont{\'e}e and Descente operators with the class ${\mathcal P}(\R^k)$ have been defined much earlier, and we follow \cite{gneiting2002} to summarize them here. The Descente and Mont{\'e}e operators are respectively defined as 
\begin{equation}\label{eq-deriv_g}
\mathcal D\varphi(t)=\left\{
\begin{array}{cc}
1, &t=0\\\\
\dfrac{\varphi'(t)}{t\varphi^{''}(0)}, & t>0,
\end{array}
\right.
\end{equation}
where $\varphi^{''}(0)$ denotes the second derivative of $\varphi$ evaluated at $t=0$, and
\begin{equation}\label{eq-integ-g}
\widetilde{\mathcal I}\varphi(t) = \int_t^\infty u\varphi(u)du \left( \int_0^\infty u\varphi(u)du\right)^{-1}.
\end{equation}
\cite{gneiting2002} proved that if $\varphi\in\mathcal P(\R^k)$, $k\geq3$, and $u\varphi(u)$ is integrable over $[0,\infty)$, then $\widetilde{\mathcal I}\varphi\in\mathcal P(\R^{k-2})$.
Invoking standard properties of Bessel functions in concert with direct inspection, \cite{gneiting2002}  proved that,
	if $\varphi\in\mathcal P(\R^k)$ and $\varphi^{''}(0)$ exists, then $\mathcal D\varphi\in\mathcal P(\R^{k+2})$. Under mild regularity conditions, the operator $\mathcal D $ and $\widetilde{\mathcal I}$ are inverse operators: 
	
	$$ \widetilde{\mathcal I} (\mathcal D \varphi ) = \mathcal D (\widetilde{\mathcal I} \varphi) = \varphi. $$

\subsection{Descente and Montée Operators on Generalized Spaces}

We start by defining the following Descente and Montée operators. The first is actually taken from 
\cite{bing-sym2019}: we define the derivate operator $D_1$ by
\begin{equation}\label{eq-descente-derv_x}
D_1\varphi(x,t) := \varphi_x(x,t)=
\dfrac{\partial \varphi}{\partial x}(x,t), \quad (x,t)\in[-1,1]\times[0,\infty).
\end{equation}
%Proposition 3.8 in \cite[]{porcu-berg} shows that $\varphi\in\mathcal P(\S^{d}\times\R^k)$ is continuously differentiable with respect to the first variable.
\\
The integral operator $I_1$ is given by
\begin{equation}\label{eq-integ-bing-x-t}
I_1\varphi(x,t) := \int_{-1}^x \varphi(u,t)du, \quad (x,t)\in[-1,1]\times[0,\infty),
\end{equation}
when $\varphi(u,t)$ is integrable over $[-1,1]$ for each $t\in[0,\infty)$.  

We define

\begin{equation}\label{eq-descente-derv_t-00}
D_2\varphi(x,t) := \left\{
\begin{array}{ll}
1, & (x,t)=(1,0)\\
\dfrac{\varphi_t(x,t)}{t\varphi_{tt}(1,0)}, & (x,t)\in[-1,1)\times(0,\infty)
\end{array}
\right. 
\end{equation}
whenever $\varphi_{tt}(1,0):=\dfrac{\partial^2 \varphi}{\partial t^2}(1,0)$ exists, and 
\begin{equation}\label{eq-integ-gne-x-t}
I_2\varphi(x,t) := \frac{\int_t^\infty v\varphi(x,v)dv}{\int_0^\infty v\varphi(1,v)dv},\quad (x,t)\in[-1,1]\times[0,\infty),
\end{equation}
when $v\varphi(1,v)$ is integrable over $[0,\infty)$ and provided the denominator is not identically equal to zero. 

The composition between the operators defined in \cite{bing-sym2019} and \cite{gneiting2002} provides a new operator, that we define here as 
	\begin{equation}\label{eq-def-I3}
	I_3\varphi(x,t) := \int_t^\infty \int_{-1}^x v\varphi(u,v)dudv,\quad  (x,t)\in[-1,1]\times[0,\infty),
\end{equation}
when $v\varphi(u,v)$ is integrable over $[-1,1]\times[0,\infty)$.%, and again provided the denominator is not identically equal to zero. 

Given $\kappa\in\Z$, we define  the operator  $I_j^\kappa$ by recurrence as:  
$$
I_j^0\varphi:=\varphi,\quad I_j^1\varphi:=I_j\varphi\quad \text{and } \quad I_j^\kappa\varphi := I_j (I_j^{\kappa-1}\varphi), \quad j=1,2,3. 
$$

\section{Dimension Walks within the Class  ${\mathcal P}(\S^d\times\R^k$)}	
This section contains our original findings. Proofs are deferred to the Appendix. 

\subsection{Descente Operators}
We start with a simple result that is an extension of Theorem 2.3 in 	\cite{beat-zucatell-2016}. In the Appendix \ref{appendix} we provide a quick sketch of the main steps.

	\begin{thm}\label{t-derv1}
	If $\varphi:[-1,1]\times[0,\infty)\to\R$  belongs to $\mathcal P(\S^d\times\R^k)$, then
$D_1\varphi$ belongs to $\mathcal P(\S^{d+2}\times\R^k)$.
\end{thm}	
%\subsubsection{Descente Operator $D_2$}	\label{s-D2}
Next result requires instead a lengthy proof and relates about the operator $D_2$. 
\begin{thm}\label{t-derv2}
%Let $d\geq2$, $k\geq1$, 
Let $d,k\in\Z^*$, 
$\varphi:[-1,1]\times[0,\infty)\to\R$  be a function in $\mathcal P(\S^d\times\R^k)$ and let $F_n$ be the $k$-Schoenberg measures associated with the $d$-Schoenberg functions of $\varphi$. If
\begin{enumerate}[(i)]
	\item $\displaystyle \int_0^\infty r^2dF_n(r)<\infty$, for all $n\in\Z$;
\label{h-int-r2-I2}
%	\item  $\dfrac{\partial^2 \varphi}{\partial t^2}(1,0)$ exists; ESTE ITEM É EQUIVALENTE A: (VER. EQ. \eqref{eq-derv-0-0})
	\item $\displaystyle 0< \dfrac{\partial^2 \varphi}{\partial t^2}(1,0)=\sum_{n=0}^{\infty} \int_0^\infty r^2dF_n(r)<\infty$,
	\label{h-ser-r2-I2} 
\end{enumerate}
then $D_2\varphi$ belongs to $\mathcal P(\S^{d+2}\times\R^k)$.
\end{thm}

\subsection{Montée Operators}

In this section we consider  functions $\varphi:[-1,1]\times[0,\infty)\to\R$ belonging to $\mathcal P(\S^d\times\R^k)$ as in \eqref{eq-pd_SdxRk} such that
\begin{equation}
\displaystyle\int_0^\infty (1/r^{2\kappa})dF_n(r)<\infty \qquad (\kappa\in\Z^*).
\end{equation}
Thus, the functions defined by
\begin{equation}\label{eq-g_n(0-k)}
g_n^\kappa(t) := \int_0^\infty \Omega_{k-2\kappa}(tr)\frac{1}{r^{2\kappa}}dF_n(r), \quad t\in[-1,1],\quad n,\kappa\in\Z,
\end{equation}
belong to the class $\mathcal P(\R^{k-2\kappa})$.

%\subsubsection{Montée Operator $I_1$} \label{s-I1}
The first finding relates to the operator $I_1$. Again, the proof is deferred to the Appendix.
\begin{thm}\label{t-I1}
%Let $d\geq3$, $k\geq1$.	If $\varphi\in\mathcal P(\S^d\times\R^k)$ and $u\mapsto\varphi(u,t)$ is integrable over $[-1,1]$ for each $t\in[0,\infty)$, then there exists a constant $C$ such that $C+I_1 \varphi\in\mathcal P(\S^{d-2}\times\R^k)$.??????
Let $k,\kappa\in\Z^*$ and $d$ be an integer such  that $d>2\kappa$. 
If $\varphi:[-1,1]\times[0,\infty)\to\R$  is a function in $\mathcal P(\S^d\times\R^k)$ such that $u\mapsto I_1^{\kappa-1}\varphi(u,t)$ is integrable over $[-1,1]$ for each $t\in[0,\infty)$.
Then, the function $I_1^\kappa\varphi$ has a representation in a form of Gegenbauer series:
\begin{equation}\label{eq-I1k}
I_1^\kappa(x,t) = \sum_{n=0}^{\infty}\tilde f_{n}^{d,\kappa}(t)C_{n}^{(d-2\kappa-1)/2}(x), \quad (x,t)\in[-1,1]\times[0,\infty),
\end{equation}
where
\begin{equation}\label{eq-g_ndk}
\tilde f_{n}^{d,\kappa}(t) := 	\left\{
\begin{array}{lll}
\displaystyle  \tau^{d,\kappa} \sum_{i=0}^{\infty}(-1)^{i}\chi_i^{n,d,\kappa}f_{i}^d(t), &n=0,1,\ldots,\kappa-1,
 \\\\
\displaystyle  \tau^{d,\kappa} f_{n-\kappa}^d(t), & n\geq\kappa.
\end{array}	
\right.
\end{equation}
The functions $f_n^d$ are the $d$-Schoenberg functions of $\varphi$ as in \eqref{eq-pd_SdxRk}, the positive constant 
 $\tau^{d,\kappa}:=(\prod_{j=1}^\kappa\delta_{(d-2j+1)/2})^{-1}$ and the coefficients
\begin{equation}\label{eq-def-chi}
\left\{	\begin{array}{ll}
\chi_i^{0,d,\kappa} &:= \displaystyle\sum_{j=1}^{\kappa-1}
(-1)^{j+1}\chi_i^{j-1,d,\kappa-1}C_{j}^{(d-2\kappa-1)/2}(1)
-
(-1)^{\kappa+1}C_{i+\kappa}^{(d-2\kappa-1)/2}(1),
\\
\\
\chi_i^{n,d,\kappa} &:= \chi_i^{n-1,d,\kappa-1}, \quad n=1,2,\ldots,\kappa-1,
\end{array}
\right.
\end{equation}
satisfy
\begin{equation}\label{eq-chi}
|\chi_i^{n,d,\kappa}|\leq \varUpsilon^{n,d,\kappa}C_i^{(d-1)/2}(1),\quad n=0,1,\ldots,\kappa-1, \quad i\in\Z,
\end{equation}
where, for each $n=0,1,\ldots,\kappa-1$ and $\kappa\in\Z^*$,  $\varUpsilon^{n,d,\kappa}$ is a positive constant that depends only on $d$.
Moreover,
$\sum_{n=0}^\infty \tilde f_{n}^{d,\kappa}(0){C_{n}^{(d-2\kappa-1)/2}(1)}	<\infty$. \\
\end{thm}

\begin{cor}\label{c-I1}
	Under the conditions of Theorem \ref{t-I1}, there exists a bounded function $H^\kappa$ on $[-1,1]\times[0,\infty)$ such that $H^\kappa+I_1^\kappa\varphi$ belongs to $\mathcal P(\S^{d-2\kappa}\times\R^{k})$.
\end{cor}

\begin{rem}\label{r-chi-positivo}
	Direct inspection shows that 
$\chi_i^{0,d,\kappa} \geq0$, for $\kappa=1,2$. Therefore,
	$\chi_i^{1,d,\kappa},\chi_i^{2,d,\kappa}\ldots,\chi_i^{\kappa-1,d,\kappa}\geq0$ for all $\kappa\geq2$ and  $i\in\Z$.
\end{rem}

\begin{rem} 
By Remark \ref{r-chi-positivo}, if $f_{2n+1}^{d}\equiv0$ for all $n$, then  $I_1^\kappa\varphi$, for $\kappa=1,2$, belongs to the class $\mathcal P(\S^{d-2\kappa}\times\R^{k})$.
Therefore, our result generalizes the corrected version  of Theorem 2.1 in \cite{bing-sym2019}.
\end{rem}

We can modify the functions $\tilde f_{n}^{d,\kappa}$, $n=0,1,\ldots,\kappa-1$, in \eqref{eq-I1k} so that the new {\it quasi} Montée operator %$I_1^\kappa\varphi$ 
belongs to $\mathcal P(\S^{d-2\kappa}\times\R^{k})$. Theorem \ref{t-I1k-AB}  sheds some light in this direction.

\begin{thm}\label{t-I1k-AB}
Let the functions $\tilde f_{n}^{d,\kappa}\in\mathcal P(\R^{k})$, $n\geq\kappa$,  and $h_{1,n}^\kappa,h_{2,n}^\kappa\in\mathcal P(\R^{k})$ be as respectively  defined at \eqref{eq-g_ndk} and \eqref{eq-h1k-h2k}.

Let $k,\kappa\in\Z^*$ and let $d$ be an integer such  that $d>2\kappa$.
	Let $\varphi:[-1,1]\times[0,\infty)\to\R$  be a function in $\mathcal P(\S^d\times\R^k)$ such that $u\mapsto I_1^{\kappa-1}\varphi(u,t)$ is integrable over $[-1,1]$ for each $t\in[0,\infty)$. If 
	\begin{enumerate}[(i)]
%		\item \label{h-int-geg-F-finita}
%	$\displaystyle \int_0^\infty C_{n}^{(d-1)/2}(1) dF_{n}(r)<\infty$;		
		\item \label{h-serie-int-med-finita}
		$\displaystyle\sum_{n=0}^{\infty} C_{n}^{(d-1)/2}(1)
		 \int_0^\infty dF_{n}(r)<\infty$;
		\item	\label{h-serie-medida-unif-limitada}
there exists constant $K>0$ such that		$\displaystyle\sum_{n=0}^{\infty} C_{n}^{(d-1)/2}(1)
		dF_{n}(r)\leq K, \,0\leq r<\infty$,	 
	\end{enumerate}	
	then there exist $2\kappa$ constants $A^n$ and $B^n$, $n=0,\ldots,\kappa-1$, such that 
%	\begin{equation}\label{eq-I1kAB}
%	I_1^{\kappa,A,B}(x,t) := \sum_{n=\kappa}^{\infty}\tilde f_{n}^{d,\kappa}(t)C_{n}^{(d-2\kappa-1)/2}(x) + (Ah_1^\kappa(t)  - Bh_2^\kappa(t))C_{\kappa-1}^{(d-2\kappa-1)/2}(x),
%	\end{equation}
\begin{eqnarray}\label{eq-I1kAB}
I_1^{\kappa,A^{0}..A^{\kappa-1},B^{0}..B^{\kappa-1}}\varphi(x,t) &:=& 
\sum_{n=0}^{\kappa-1}
(A^nh_{1,n}^\kappa(t)  - B^nh_{2,n}^\kappa(t))C_{n}^{(d-2\kappa-1)/2}(x) 
\nonumber 
\\ 
&+& \sum_{n=\kappa}^{\infty}\tilde f_{n}^{d,\kappa}(t)C_{n}^{(d-2\kappa-1)/2}(x) ,
\end{eqnarray}
belongs to $\mathcal P(\S^{d-2\kappa}\times\R^{k})$.

\end{thm}

\begin{rem} For any $A^n\geq0$, $n=1,\ldots,\kappa-1$ the function $I_1^{\kappa,0A^1..A^{\kappa-1},0..0}\varphi$ belongs to $\mathcal P(\S^{d-2\kappa}\times\R^{k})$. This also can be seen as a generalization of the correction of Theorem 2.1 in \cite{bing-sym2019} (to appear).
\end{rem}

Next results is related to the operator $I_2$. 
\begin{thm} \label{t-I2_0-k}
%	Let $d\in\Z^*$, $\kappa\in\Z^*$ and $k$ be an integer greater than $3+2\kappa$. 
Let  $d,\kappa\in\Z^*$ and $k$ be an integer such  that $k>2\kappa$.
If $\varphi:[-1,1]\times[0,\infty)\to\R$  is a function in $\mathcal P(\S^d\times\R^k)$ such that  
\begin{enumerate}[(i)]
\item $g_n^\nu(0)=\displaystyle\int_0^\infty (1/r^{2\nu})dF_n(r)<\infty$, for all $n\in\Z$ and $\nu\in\{1,2,\ldots,\kappa\}$. % \rosa{tem que ajustar onde corre $\nu$}
	\label{h-int-r2-I2_0-k}
\item $0\neq\sum_{n=0}^{\infty}g_n^{(\nu)}(0){C_{n}^{(d-1)/2}(1)}
%= \displaystyle \sum_{n=0}^\infty\int_0^\infty (1/r^{2\nu})dF_n(r)
<\infty$, for $\nu\in\{1,2,\ldots, \kappa\}$,
	\label{h-serie-I2_0-k}	  
\end{enumerate}
then the function $I_2^\kappa\varphi$ has a representation in Gegenbauer series in the form
	\begin{eqnarray}\label{eq-I2_0-k}
	I_2^\kappa\varphi(x,t) = 
	\frac{1}{ \sum_{n=0}^{\infty}g_n^\kappa(0){C_{n}^{(d-1)/2}(1)}} \sum_{n=0}^{\infty}g_n^\kappa(t) C_{n}^{(d-1)/2}(x).
	\end{eqnarray}
The functions $g_n^\kappa$ are defined in \eqref{eq-g_n(0-k)}
%	\begin{equation}\label{eq-g_n(0-k)}
%	g_n^\kappa(t) := \int_0^\infty \Omega_{k-2\kappa}(tr)\frac{1}{r^{2\kappa}}dF_n(r),
%	\end{equation}
and $F_n$ are the $k$-Schoenberg measures   of the $d$-Schoenberg functions of $\varphi$. 

Moreover,  $I_2^\kappa \varphi$ belongs to $\mathcal P(\S^{d}\times\R^{k-2\kappa})$.
\end{thm}

We finish this part with the Mont{\'e}e operator $I_3$.

\begin{thm}\label{t-I3k+H}
%	Let $d\geq2$, $\kappa\in\Z^*$ and $k$ be an integer greater than $3+2\kappa$. 
Let $\kappa\in\Z^*$, $d$ and $k$ be integers such  that $d,k>2\kappa$. %and  $k\geq3+2\kappa$.
If $\varphi:[-1,1]\times[0,\infty)\to\R$  is a function in $\mathcal P(\S^d\times\R^k)$ such that 
\begin{enumerate}[(i)]	
\item $g_n^{\nu}(0)=\displaystyle\int_0^\infty \frac{1}{r^{2\nu}}dF_n(r) <\infty$ for all $n\in\Z$ and $\nu\in\{1,2,\ldots,\kappa\}$;
\label{h-int_r2-I3k}
\item $\displaystyle \sum_{n=0}^\infty g_{n}^\nu(0)=\sum_{n=0}^\infty\int_0^\infty \frac1{r^{2\nu}}dF_{n}(r)<\infty$, for and $\nu\in\{1,2,\ldots,\kappa\}$,
\label{h-serie-gn-finit-I3k}
\item $\displaystyle {\sum_{n=0}^\infty g_{n}^\nu(0)C_{n+\nu}^{(d-2\nu-1)/2}(1)}
%=\sum_{n=0}^\infty\int_0^\infty \frac1{r^{2\nu}}dF_{n}(r)
<\infty$, for and $\nu\in\{1,2,\ldots,\kappa\}$,
\label{h-serien-finit-I3k}
	\label{h-serie-nozero-I3k}
\end{enumerate}
then the function $I_3^\kappa\varphi$ has a representation in Gegenbauer series in the form
	\begin{equation}\label{eq-I3k}
	I_3^\kappa\varphi(x,t) = 
	%\frac1{\gamma^{(d-2\kappa-1)/2,\kappa}}
	\sum_{n=0}^{\infty} h_n^\kappa(t)C_{n}^{(d-2\kappa-1)/2}(x), \quad (x,t)\in[-1,1]\times[0,\infty),
	\end{equation}
where
	\begin{eqnarray}\label{eq-def-h_n(k)}
	h_n^\kappa(t):= \left\{ 
	\begin{array}{lll} 
	\displaystyle	
	%\frac{1}{2\sum_{j=0}^{\infty} g_{2j}^\kappa(0)} 
	\gamma^{d,k,\kappa}
	\sum_{i=0}^{\infty}(-1)^{i}\chi_i^{n,d,\kappa}g_{i}^{\kappa}(t), & n=0,1,\ldots,\kappa-1, \\ \\
	\displaystyle	 
	%\frac{1}{2\sum_{j=0}^{\infty} g_{2j}^\kappa(0)}
	\gamma^{d,k,\kappa}
	 g_{n-\kappa}^\kappa(t), & n\geq\kappa, \\ \\
	\end{array}	
	\right.
	\end{eqnarray}
with  %the positive constants
$\gamma^{d,k,\kappa}:=\prod_{j=1}^\kappa\frac{k-2j}{\delta_{(d-2j+1)/2}}>0$ and
{$%0\neq\gamma^{(d-2\kappa-1)/2,\kappa}:= 
	\sum_{n=0}^{\infty} h_{n}^\kappa(0)C_{n}^{(d-2\kappa-1)/2}(1)<\infty$}.
The functions $g_n^\kappa$ are defined in \eqref{eq-g_n(0-k)} and belong to the class $\mathcal P(\R^k)$ and $\chi_i^{n,d,\kappa}$ are given in \eqref{eq-def-chi}.
\end{thm}

\begin{rem} 
	By Remark \ref{r-chi-positivo}, if %\ana{$\gamma^{(d-2\kappa-1)/2,\kappa}>0$????} and 
	 $g_{2n+1}^\kappa\equiv0$ for all $n$, then  $I_3^\kappa\varphi$, for $\kappa=1,2$, belongs to
	the class $\mathcal P(\S^{d-2\kappa}\times\R^{k-2\kappa})$.
	%	If $g_{2n+1}^\kappa\equiv0$ for all $n\in\Z$  (this is a similar hypothesis of the correction of \cite{bing-sym2019}, to appear) then, by \eqref{eq-I3kAB} and Theorem 3.3 in \cite{porcu-berg}, we conclude that the Montée operator $I_3^\kappa\varphi$ belongs to the class $\mathcal P(\S^{d-2\kappa}\times\R^{k-2\kappa})$.
\end{rem}

\begin{cor} \label{c-tI3k+H}
	Under the conditions of Theorem \ref{t-I3k+H}, % if \ana{$\gamma^{(d-2\kappa-1)/2,\kappa}>0$???}, 	then	
	there exists a bounded function $H^\kappa$ on $[-1,1]\times[0,\infty)$ such that $H^\kappa+I_3^\kappa\varphi$ belongs to $\mathcal P(\S^{d-2\kappa}\times\R^{k-2\kappa})$.
\end{cor}

As previously mentioned, we can replace the functions $h_{n}^\kappa$, $n=0,1,\ldots,\kappa-1$, %in \eqref{eq-I3k} 
with others  such that the new {\em quasi} Montée operator belongs to $\mathcal P(\S^{d-2\kappa}\times\R^{k-2\kappa})$. Theorem \ref{t-I3k-pd-A-B}  provides a construction in this sense.

\begin{thm} \label{t-I3k-pd-A-B}
%Let $d\geq2$, $\kappa\in\Z^*$ and $k$ be an integer greater than $3+2\kappa$. 
Let $\kappa\in\Z^*$, $d$ and $k$ be integers such  that $d,k>2\kappa$. %and  $k\geq3+2\kappa$.
Let $\varphi:[-1,1]\times[0,\infty)\to\R$  be a function belongs to the class $\mathcal P(\S^d\times\R^k)$  satisfying the hypotheses of Theorem \ref{t-I3k+H}. 
If additionally, %\ana{$\gamma^{(d-2\kappa-1)/2,\kappa}>0$?? }and
 the $k$-Schoenberg measures $F_n$ of the $d$-Schoenberg functions of $\varphi$ satisfy 
\begin{enumerate}[(i)]
\item $\displaystyle \sum_{n=0}^\infty g_{n}^\nu(0)C_n^{(d-1)/2}(1)=\sum_{n=0}^\infty C_n^{(d-1)/2}(1)\int_0^\infty \frac1{r^{2\nu}}dF_{n}(r)<\infty$, for and $\nu\in\{1,2,\ldots,\kappa\}$,
\label{h-serie-gn-Cn-finit-I3k}	
%	\item $\displaystyle\sum_{n=0}^{\infty}\frac{1}{r^{2\kappa}}dF_{n}(r)<\infty;$
%	\label{h-serie-dF_n}
\item there exists constant $K>0$ such that $\displaystyle\sum_{n=0}^{\infty}C_n^{(d-1)/2}(1)dF_{n}(r)\leq K$, $0\leq r<\infty$;
	\label{eq-cond-medida_F_n}
\end{enumerate}
then there exist exist $2\kappa$ constants $A^n$ and $B^n$, $n=0,\ldots,\kappa-1$, such that 
\begin{eqnarray}\label{eq-I3kAB}
	I_3^{\kappa,A^{0}..A^{\kappa-1},B^{0}..B^{\kappa-1}}\varphi(x,t) &:=& 
	\sum_{n=0}^{\kappa-1}
	(A^n\tilde h_{1,n}^\kappa(t)  - B^n\tilde h_{2,n}^\kappa(t))C_{n}^{(d-2\kappa-1)/2}(x) 
	\nonumber 
	\\ 
	&+& %\frac1{\gamma^{(d-2\kappa-1)/2,\kappa}}
	 \sum_{n=\kappa}^{\infty}\tilde h_{n}^{\kappa}(t)C_{n}^{(d-2\kappa-1)/2}(x).
\end{eqnarray}
The functions $h_n^\kappa\in\mathcal P(\R^k)$, $n\geq\kappa$ and $\tilde h_{1,n}^\kappa,\tilde h_{2,n}^\kappa\in\mathcal P(\R^k)$ are defined respectively in \eqref{eq-def-h_n(k)} and \eqref{eq-til-h1k-h2k}.

\end{thm}

\begin{rem} For any $A^n\geq0$, $n=1,\ldots,\kappa-1$ the function $I_3^{\kappa,0A^1..A^{\kappa-1},0..0}\varphi$ belongs to $\mathcal P(\S^{d-2\kappa}\times\R^{k-2\kappa})$. 
\end{rem}

\section*{Declarations}

{\bf Funding.}
Ana Paula Peron was partially supported by Funda\c c\~ao de Amparo \`a Pesquisa do Estado de S\~ao Paulo - FAPESP, Brazil, \# 2021/04269-0.
\\
\\
{\bf Conflict of interest.}
There were no competing interests to declare which arose during the preparation or publication process of this article.

\appendix

\section*{Appendix A}\label{appendix}

\begin{proof}[{\bf Proof of Theorem \ref{t-derv1}}]
		Since $\varphi$ belongs to $\mathcal P(\S^d\times\R^k)$, then $\varphi$ is continuously differentiable with respect to the first variable (see \cite[Propositon 3.8]{porcu-berg}) and has a Gegenbauer expansion as \eqref{eq-los_grandes_BP}. Also $\varphi_x$ has a Gegenbauer expansion in the form
		\begin{equation*}
			\varphi_x(x,t) \sim \sum_{n=0}^{\infty}\tilde f_{n}^{d+1}(t)C_{n}^{(d+1)/2}(x).
		\end{equation*} 
		Using \eqref{eq-phi__phi_x}-\eqref{eq-coef-phi__phi_x}, remainder of the proof follows as in \cite[Theorem 2.3]{beat-zucatell-2016}.
\end{proof}

\begin{proof}[{\bf Proof of Theorem \ref{t-derv2}}]
Let $\varphi$ be a function as in \eqref{eq-pd_SdxRk}. By \eqref{eq-derv-Omega},
$$
\dfrac{df_{n}^d}{dt}(t) = \int_0^\infty -\frac1k tr^2 \Omega_{k+2}(tr)dF_n(r).
$$

Deriving term by term, we obtain
\begin{eqnarray}\label{eq-derv_x_t}
\dfrac{\partial \varphi}{\partial t}(x,t)  
&=& \sum_{n=0}^{\infty}\dfrac{df_{n}^d}{dt}(t)C_{n}^{(d-1)/2}(x)\nonumber\\
%&=& \sum_{n=0}^{\infty} \left(\int_0^\infty -\frac1k tr^2 \Omega_{k+2}(tr)dF_n(r)\right) C_{n}^{(d-1)/2}(\cos x)\\
&=& -\frac1k t \sum_{n=0}^{\infty} \left(\int_0^\infty  \Omega_{k+2}(tr)r^2 dF_n(r)\right) C_{n}^{(d-1)/2}(x).
\end{eqnarray}

By Lemma 3 in \cite{gneiting99}, we have  
$$
\dfrac{d^2f_{n}^d}{dt^2}(0)= -\frac1k  \int_0^\infty  r^2 dF_n(r).
$$

Thus,
\begin{eqnarray}\label{eq-derv_x_0}
\dfrac{\partial^2 \varphi}{\partial t^2}(x,0)  &=& 
-\frac1k \sum_{n=0}^{\infty} \left(\int_0^\infty r^2dF_n(r)\right) C_{n}^{(d-1)/2}(x), \quad x\in[-1,1].
\end{eqnarray}

In particular, 
\begin{eqnarray}\label{eq-derv-0-0}
\dfrac{\partial^2 \varphi}{\partial t^2}(1,0)  &=& 
-\frac1k \sum_{n=0}^{\infty} \left(\int_0^\infty r^2dF_n(r)\right) C_{n}^{(d-1)/2}(1).
\end{eqnarray}

Thus, by \eqref{eq-derv_x_t} and \eqref{eq-derv-0-0}, for %$x\in[-1,1]\setminus\{0\}$ 
$x\in[-1,1)$ 
and $t>0$ we have
\begin{eqnarray}\label{eq-expressao-fora-0}
D_2\varphi(x,t) = \dfrac{\varphi_t(x,t)}{t\varphi_{tt}(1,0)} % \nonumber\\&=&  
= \frac1{\sum_{n=0}^{\infty} \int_0^\infty r^2dF_n(r)} \sum_{n=0}^{\infty} g_n^d(t) C_{n}^{(d-1)/2}(x), 
\end{eqnarray}
where the functions $ g_n^d(t):[0,\infty)\to\R$ are defined by
$$
g_n^d(t) := \int_0^\infty  \Omega_{k+2}(tr)r^2 dF_n(r).
$$
We invoke Hypothesis \eqref{h-int-r2-I2} to imply that $g_n^d\in\mathcal P(\R^{k+2})$. Thus, the series in the Equation \eqref{eq-expressao-fora-0} converges absolutely and uniformly on $[-1,1]\times[0,\infty)$. Hence, letting $x=1$ and $t=0$ %$x=t=0$ 
in the expression  of the series in  \eqref{eq-expressao-fora-0} provides
\begin{eqnarray*}
\frac1{\sum_{n=0}^{\infty} \int_0^\infty r^2dF_n(r)} \sum_{n=0}^{\infty} g_n^d(0) C_{n}^{(d-1)/2}(1) = 1 = D_2\varphi(1,0).
\end{eqnarray*}

Therefore, $D_2\varphi$ is a  continuous function on $[-1,1]\times[0,\infty)$ having a representation series as in \eqref{eq-pd_SdxRk} with $d$-Schoenberg functions $g_n^d\in\mathcal P(\R^{k+2})$. Since, by \eqref{h-ser-r2-I2}, 
\begin{equation*}\label{eq-serie-int-conv}
\sum_{n=0}^{\infty} g_n^d(0)= \sum_{n=0}^{\infty} \int_0^\infty r^2dF_n(r)<\infty,
\end{equation*}
we can conclude that $D_2\varphi$ belongs to $\mathcal P(S^d\times\R^{k+2})$.	

\end{proof}

\begin{proof}[{\bf Proof of Theorem \ref{t-I1}}] 
	We prove the statement by induction on $\kappa \in \Z^*$.\\
{\bf Step $\kappa=1$:} We have
\begin{eqnarray*}
I_1^1\varphi(x,t) &=&  I_1\varphi(x,t) =  \int_{-1}^x \varphi(u,t)du.
\end{eqnarray*}
By \eqref{eq-pd_SdxRk} and \eqref{eq-integ-Gegenb}, integrating term by term, we obtain
\begin{eqnarray*}
I_1^1\varphi(x,t)  &=& 
\sum_{n=0}^{\infty}f_{n}^d(t) \frac1{\delta_{(d-1)/2}}\left(C_{n+1}^{(d-3)/2}(x)-C_{n+1}^{(d-3)/2}(-1)\right).
\end{eqnarray*} 
Since   $C_{n+1}^{(d-3)/2}(-1)=(-1)^{n+1}C_{n+1}^{(d-3)/2}(1)$, we have
\begin{eqnarray*}
	I_1^1\varphi(x,t)  &=& 
\sum_{n=0}^{\infty}\tilde f_{n}^{d,1}(t)C_{n}^{(d-3)/2}(x)
\end{eqnarray*}
where
\begin{eqnarray*}
\tilde f_{n}^{d,1}(t):= 
\left\{
\begin{array}{ll}
\displaystyle  \frac1{\delta_{(d-1)/2}}\sum_{i=0}^{\infty}(-1)^{i}\chi_i^{n,d,1}f_{i}^d(t) , &n=0 \\\\
\displaystyle  \frac1{\delta_{(d-1)/2}}f_{n-1}^d(t), & n\geq1,
\end{array}	
\right.
\end{eqnarray*}	

where $\chi_i^{0,d,1} = C_{i+1}^{(d-3)/2}(1)$ and, by \eqref{eq-razao-gegenb},
\begin{eqnarray*}
0\leq \chi_i^{0,d,1} 
=
	\left|C_{i}^{(d-1)/2}(1)
	%\frac1{\delta_{(d-1)/2}}
	\frac{C_{i+1}^{(d-3)/2}(1)}{C_{i}^{(d-1)/2}(1)}\right| 
	\leq 
	\underbrace{\varrho_{(d-1)/2,1}}_{\varUpsilon^{0,d,1}}
	C_{i}^{(d-1)/2}(1),
\end{eqnarray*}
which implies
$$
\left|(-1)^{i}\chi_i^{0,d,1}f_{i}^d(t)\right| \leq \varUpsilon^{0,d,1}f_i(0)C_{i}^{(d-1)/2}(1)
$$
and by \eqref{eq-pd_SdxRk}, the series in the definition of $\tilde f_{n}^{d,1}$ is uniformly convergent on $[0,\infty)$. 
Again by \eqref{eq-razao-gegenb}, for $n\geq1$,
\begin{eqnarray*}
\left|
%\frac1{\delta_{(d-1)/2}}
\tilde f_{n}^{d,1}(0)C_{n}^{(d-3)/2}(1)\right| 
%= 
%\frac1{\delta_{(d-1)/2}}
%\left| f_{n-1}^d(0)C_{n-1}^{(d-1)/2}(1)
%\frac{C_{n}^{(d-3)/2}(1)}{C_{n-1}^{(d-1)/2}(1)}\right| 
\leq 
\frac1{\delta_{(d-1)/2}}\varUpsilon^{0,d,1}f_{n-1}^d(0)C_{n-1}^{(d-1)/2}(1).
\end{eqnarray*}
Thus, 
$
\sum_{n=0}^{\infty}\tilde f_{n}^{d,1}(0)C_{n}^{(d-3)/2}(1)<\infty$.\\

%%%%%%%%%%%%%%%%%%% PASSO k=2
%

{\bf Step $\kappa=2$:}  By algebraic manipulation we have

\begin{eqnarray*}
	I_1^2\varphi(x,t)  &=& 
	\sum_{n=0}^{\infty}\tilde f_{n}^{d,2}(t)C_{n}^{(d-5)/2}(x)
\end{eqnarray*}
where
\begin{equation*}
	\tilde f_{n}^{d,2}(t) := 	\left\{
	\begin{array}{lll}
\tau^{d,2}
	\sum_{i=0}^{\infty}(-1)^{i}\chi_i^{n,d,2}f_{i}^d(t)
		, & n=0, 1, \quad \\\\
		\displaystyle  \tau^{d,2} f_{n-\kappa}^d(t), & n\geq2,
	\end{array}	
	\right.
\end{equation*}
with $\tau^{d,2}=(\delta_{(d-1)/2}\delta_{(d-3)/2})^{-1}$ and
$$
\chi_i^{0,d,2}:= C_{i+1}^{(d-3)/2}(1) C_{1}^{(d-5)/2}(1) -C_{i+2}^{(d-5)/2}(1),
$$
$$
\chi_i^{1,d,2}:=C_{i+1}^{(d-3)/2}(1) = \chi_i^{0,d,1}.
$$

It is clear that $0\leq \chi_i^{1,d,2}\leq \varUpsilon^{1,d,2}C_i^{(d-1)/2}(1)$, with $\varUpsilon^{1,d,2}%=\varrho_{(d-1)/2,1}
=\varUpsilon^{0,d,1}$. It is not difficult to see that
$$
\chi_i^{0,d,2} = \frac{(d-5)\Gamma(i+d-1)(i+1)}{\Gamma(d-3)(i+d-3)\Gamma(i+3)} \geq0
$$
and, by \eqref{eq-razao-gegenb},
\begin{eqnarray*}
|\chi_i^{0,d,2}| 
&\leq&
\left(
\left|\frac{C_{i+1}^{(d-3)/2}(1) }{C_i^{(d-1)/2}(1)} C_{1}^{(d-5)/2}(1)\right| + \left|\frac{C_{i+2}^{(d-5)/2}(1)}{C_i^{(d-1)/2}(1)}\right|
\right)
C_i^{(d-1)/2}(1)
\\
&\leq&
\left(\varrho_{(d-1)/2,1}\frac{\Gamma(1+d-5)}{\Gamma(2)\Gamma(d-5)} + \varrho_{(d-1)/2,2} 
\right)C_i^{(d-1)/2}(1)
\\
&=&
\underbrace{\left(\varrho_{(d-1)/2,1}(d-5) + \varrho_{(d-1)/2,2}\right) }_{\varUpsilon^{0,d,2}}
C_i^{(d-1)/2}(1)
\end{eqnarray*}

By the same argument of the step $\kappa=1$ we can conclude that the series in the definition of $\tilde f_{n}^{d,2}$ for $n=0,1$  are uniformly convergent on $[0,\infty)$  and also $\sum_{n=0}^{\infty}\tilde f_{n}^{d,2}(0)C_{n}^{(d-5)/2}(1)<\infty$.

\vspace{2mm}
	
\noindent {\bf Induction step:} let us assume that the expression in \eqref{eq-I1k} of $I_1^\kappa\varphi$ holds up to $\kappa$, and let us prove it holds for $I_1^{\kappa+1}\varphi$.
	We have
$$
I_1^{\kappa+1}\varphi(x,t) = I_1(I_1^\kappa\varphi)(x,t) = 
\int_{-1}^x  I_1^\kappa\varphi(u,t)du.
$$
	
Using the induction hypothesis and integrating term by term,
% using \eqref{eq-g_ndk} and \eqref{eq-integ-Gegenb}, 
for $(x,t)\in[-1,1]\times[0,\infty)$, we obtain
\begin{eqnarray*}
&&I_1^{\kappa+1}\varphi(x,t) 
=
\sum_{n=0}^{\kappa-1}\tilde f_{n}^{d,\kappa}(t)\int_{-1}^x C_{n}^{(d-2\kappa-1)/2}(u)du
+
\sum_{n=\kappa}^{\infty}\tilde f_{n}^{d,\kappa}(t)\int_{-1}^x C_{n}^{(d-2\kappa-1)/2}(u)du
\\
&=&
\sum_{n=0}^{\kappa-1}
\tau^{d,\kappa}\sum_{i=0}^{\infty}(-1)^{i}\chi_i^{n,d,\kappa}f_{i}^d(t)
\frac1{\delta_{(d-2\kappa-1)/2}}
\left(C_{n+1}^{(d-2\kappa-3)/2}(x)-C_{n+1}^{(d-2\kappa-3)/2}(-1)\right)
\\
&&
+
\sum_{n=\kappa}^{\infty}
%\frac1{\prod_{j=1}^\kappa\delta_{(d-2j+1)/2}}
\tau^{d,\kappa}
f_{n-\kappa}^d(t)
\frac1{\delta_{(d-2\kappa-1)/2}}
\left(C_{n+1}^{(d-2\kappa-3)/2}(x)-C_{n+1}^{(d-2\kappa-3)/2}(-1)\right).
\end{eqnarray*}

Thus,
\begin{eqnarray*}
%	&&
	I_1^{\kappa+1}\varphi(x,t) 
%	=
%\\
&=&	
	\tau^{d,\kappa+1}
	\sum_{n=0}^{\kappa-1}
	\sum_{i=0}^{\infty}(-1)^{i}\chi_i^{n,d,\kappa}f_{i}^d(t)
	\left(C_{n+1}^{(d-2\kappa-3)/2}(x)-C_{n+1}^{(d-2\kappa-3)/2}(-1)\right)
	\\
	&&
	+
	\tau^{d,\kappa+1}
	\sum_{n=\kappa}^{\infty}
	f_{n-\kappa}^d(t)
	\left(C_{n+1}^{(d-2\kappa-3)/2}(x)-C_{n+1}^{(d-2\kappa-3)/2}(-1)\right).
\end{eqnarray*}

After some algebraic manipulation,
$$
I_1^{\kappa+1}\varphi(x,t) 
=
\sum_{n=0}^{\infty}\tilde f_{n}^{d,\kappa+1}(t)C_{n}^{(d-2(\kappa+1)-1)/2}(x),
$$
where
\begin{equation*}
\tilde f_{n}^{d,\kappa+1}(t) = 	\left\{
\begin{array}{lll}
\displaystyle  
\tau^{d,\kappa+1}
\sum_{i=0}^{\infty}(-1)^{i}\chi_i^{n,d,\kappa+1}f_{i}^d(t), 
&n=0,1,\ldots,\kappa,  \\\\
\displaystyle 
\tau^{d,\kappa+1}
f_{n-(\kappa+1)}^d(t), & n\geq\kappa+1,
\end{array}	
\right.
\end{equation*}	
with
$$
\chi_i^{0,d,\kappa+1} := \sum_{j=1}^{\kappa}
(-1)^{j+1}\chi_i^{j-1,d,\kappa}C_{j}^{(d-2(\kappa+1)-1)/2}(1)
-
(-1)^{\kappa+1}C_{i+\kappa+1}^{(d-2(\kappa+1)-1)/2}(1)
$$
and
$$
\chi_i^{n,d,\kappa+1} := \chi_i^{n-1,d,\kappa}, \quad n=1,2,\ldots,\kappa.
$$

It is clear that $|\chi_i^{n,d,\kappa+1}|\leq \varUpsilon^{n,d,\kappa+1}C_i^{(d-1)/2}(1)$, with $\varUpsilon^{n,d,\kappa+1}=\varUpsilon^{n-1,d,\kappa}$. Now, by \eqref{eq-razao-gegenb},
\begin{eqnarray*}
|\chi_i^{0,d,\kappa+1}|
&\leq& 
\left( 
\sum_{j=1}^{\kappa}
\frac{|\chi_i^{j-1,d,\kappa}|}{C_{i}^{(d-1)/2}(1)}C_{j}^{(d-2\kappa-3)/2}(1)
+
\frac{C_{i+\kappa+1}^{(d-2\kappa-3)/2}(1)}{C_{i}^{(d-1)/2}(1)}
\right)C_{i}^{(d-1)/2}(1).
\end{eqnarray*}

By induction hypothesis \eqref{eq-chi} and by \eqref{eq-razao-gegenb} we obtain

\begin{eqnarray*}
|\chi_i^{0,d,\kappa+1}|
	&\leq& 
\underbrace{\left( 
\sum_{j=1}^{\kappa}
\varUpsilon^{j-1,d,\kappa}(d-2\kappa-3)
	+
\varrho_{(d-1)/2,2(\kappa+1)}	
	\right)}_{\varUpsilon^{0,d,\kappa+1}}
C_{i}^{(d-1)/2}(1).
\end{eqnarray*}

The convergence of the series  in the definition of $\tilde f_{n}^{d,\kappa+1}$ for $n=0,1,\ldots, \kappa$ and $\sum_{n=0}^{\infty}\tilde f_{n}^{d,\kappa+1}(0)C_{n}^{(d-2(\kappa+1)-1)/2}(1)$ follow as in the previous steps.

\end{proof}	

\begin{proof}[{\bf Proof of Corollary \ref{c-I1}}] 

%Given $\kappa\in\Z^*$, by \eqref{eq-pd_SdxRk}, it is clear that $\tilde f_{n}^{d,\kappa}\in\mathcal P(\R^{k})$, for all $n\geq\kappa$. 

Note that, for $n=0,1,\ldots,\kappa-1$, we can rewrite $\tilde f_{n}^{d,\kappa}$  as 
$$
\tilde f_{n}^{d,\kappa}(t) = h_{1,n}^\kappa(t) - h_{2,n}^\kappa(t),
$$
where 
\begin{eqnarray}\label{eq-h1k-h2k}
	&&
	h_{1,n}^\kappa(t) := \tau^{d,\kappa} \sum_{i=0}^{\infty}\chi_{2i}^{n,d,\kappa}f_{2i}^d(t),\qquad \text{and} 
	\\
	&&
	h_{2,n}^\kappa(t) := \tau^{d,\kappa} \sum_{i=0}^{\infty}\chi_{2i+1}^{n,d,\kappa}f_{2i+1}^d(t)
\end{eqnarray}

Define the function $H^\kappa$ on $[-1,1]\times[0,\infty)$ by 
$$
H^\kappa(x,t):=  \sum_{n=0}^{\kappa-1}h_{2,n}^{\kappa}(t)C_{n}^{(d-2\kappa-1)/2}(x) - h_{1,0}^{\kappa}(t)C_{0}^{(d-2\kappa-1)/2}(x) 
$$
which is bounded on $[-1,1]\times[0,\infty)$  because, by \eqref{eq-Gegenb-ltda},  \eqref{eq-chi} and \eqref{eq-pd_SdxRk},
\begin{eqnarray*}
	|H^\kappa(x,t)| 
%	&\leq&
%\tau^{d,\kappa} \sum_{n=0}^{\kappa-1} \sum_{i=0}^{\infty}|\chi_{2i+1}^{n,d,\kappa}|\, |f_{2i+1}^d(t)| C_{n}^{(d-2\kappa-1)/2}(1)
%+
%\tau^{d,\kappa} \sum_{i=0}^{\infty}|\chi_{2i}^{0,d,\kappa}|\,|f_{2i}^d(t)|C_{0}^{(d-2\kappa-1)/2}(1)
%\\
&\leq&
\tau^{d,\kappa} \sum_{n=0}^{\kappa-1}
\varUpsilon^{n,d,\kappa}
\left( \sum_{i=0}^{\infty} f_{2i+1}^d(0) C_{2i+1}^{(d-1)/2}(1)\right)C_{n}^{(d-1)/2}(1)
\\
&&
+
\tau^{d,\kappa}\varUpsilon^{0,d,\kappa} 
\left(\sum_{i=0}^{\infty} f_{2i}^d(0)C_{2i}^{(d-1)/2}(1)\right)
C_{0}^{(d-2\kappa-1)/2}(1)
<\infty, 
\end{eqnarray*}
for all $(x,t)\in[-1,1]\times[0,\infty)$.

By Remark \ref{r-chi-positivo}, it is clear that $h_{1,n}^\kappa\in\mathcal P(\R^{k})$, $n=1,2,\ldots,\kappa-1$, and also $\tilde f_{n}^{d,\kappa}\in\mathcal P(\R^{k})$ for $n\geq\kappa$. Therefore,
$$
H^\kappa(x,t)+I_1^\kappa\varphi(x,t) =  
\sum_{n=1}^{\kappa-1}h_{1,n}^\kappa(t) C_{n}^{(d-2\kappa-1)/2}(x) 
+
\tau^{d,\kappa}\sum_{n=\kappa}^{\infty}\tilde f_{n}^{d,\kappa}(t)C_{n}^{(d-2\kappa-1)/2}(x)
$$
has an expansion uniformly convergent as \eqref{eq-pd_SdxRk}
due to Theorem \ref{t-I1}.
By Theorem 3.3 of \cite{porcu-berg} (see \eqref{eq-pd_SdxRk}), we can conclude that the function $H^\kappa+I_1^\kappa\varphi$ belongs to the class $\mathcal P(\S^{d-2\kappa}\times\R^{k})$.

\end{proof}	

{We observe that the function $H^\kappa$ is not unique and that the construction presented allows us to highlight the properties of the coefficient functions and consider the maximum of the non-zero $d$-Schoenberg functions of $\varphi$.}

\begin{proof}[{\bf Proof of Theorem \ref{t-I1k-AB}}]
By \eqref{eq-pd-Rk}, for any  constants $A$ and $B$,
\begin{equation}\label{eq-serieh1-h2-I1}
Ah_{1,n}^\kappa(t)  - Bh_{2,n}^\kappa(t) 
	=
\tau^{d,\kappa}
	\sum_{i=0}^{\infty} \int_0^\infty  \Omega_k(tr)
	\left(A\chi_{2i}^{n,d,\kappa}dF_{2i}(r) - B\chi_{2i+1}^{n,d,\kappa}dF_{2i+1}(r)\right).
\end{equation}	

Since, by \eqref{eq-chi},
$$
\int_0^\infty A|\chi_{i}^{n,d,\kappa}|%C_{i}^{(d-1)/2}(1)
dF_{i}(r) 
\leq
A\varUpsilon^{n,d,\kappa}C_{i}^{(d-1)/2}(1)
\int_0^\infty dF_{i}(r) ,
$$
we have
$$
\int_0^\infty\left(A\chi_{2i}^{n,d,\kappa}dF_{2i}(r) - B\chi_{2i+1}^{n,d,\kappa}dF_{2i+1}(r)\right)<\infty.
$$
By \eqref{eq-chi} and \eqref{h-serie-int-med-finita} the series  in \eqref{eq-serieh1-h2-I1} converges absolutely and uniformly on $[0,\infty)$. 

Thus,
\begin{eqnarray*}
	Ah_{1,n}^\kappa(t)  - Bh_{2,n}^\kappa(t) 
	&=&
	\tau^{d,\kappa}
	 \int_0^\infty  \Omega_k(tr)
	d\left(	
	\sum_{i=0}^{\infty} 
		A\chi_{2i}^{n,d,\kappa}F_{2i}(r) - B\chi_{2i+1}^{n,d,\kappa}F_{2i+1}(r)\right).
\end{eqnarray*}	

By  \eqref{eq-chi} and \eqref{h-serie-medida-unif-limitada},  the series $\sum_{i=0}^{\infty} 
\chi_{2i}^{n,d,\kappa}F_{2i}$ and $\sum_{i=0}^{\infty}  \chi_{2i+1}^{n,d,\kappa}F_{2i+1}$ are uniformly bounded on $[0,\infty)$. 
Then we can choose $A^n,B^n$ such that the series
 $\sum_{i=0}^{\infty} 
 A^n\chi_{2i}^{n,d,\kappa}F_{2i} - B^n\chi_{2i+1}^{n,d,\kappa}F_{2i+1}$ is non negative, which allows us conclude that  $A^nh_{1,n}^\kappa-B^nh_{2,n}^\kappa\in\mathcal P(\R^{k})$.
The convergence uniform of the series \eqref{eq-I1kAB} follows by Theorem \ref{t-I1} and the result by Theorem 3.3 of \cite{porcu-berg} (see \eqref{eq-pd_SdxRk}).
\end{proof}

\begin{proof}[{\bf Proof of Theorem \ref{t-I2_0-k}}]
	We will prove \eqref{eq-I2_0-k} by mathematical induction on $\kappa$.
\\
%	{\bf Step $\kappa=0$:} Follows from \cite[Theorem 3.3]{porcu-berg} (see \eqref{eq-pd_SdxRk}).
%\\
%
	{\bf Step $\kappa=1$:}  We have
	\begin{eqnarray*}
		I_2\varphi(x,t) =
		\frac1{\int_0^\infty v\varphi(1,v)dv} \int_t^\infty v\varphi(x,v)dv.
	\end{eqnarray*}	
	
By \eqref{eq-pd_SdxRk}, integrating term by term, we obtain
\begin{eqnarray*}
	\int_t^\infty v\varphi(x,v)dv &=& 
%	\int_t^\infty v\sum_{n=0}^{\infty}f_{n}^d(v)C_{n}^{(d-1)/2}(x)dv\\
%	&=&
	 \sum_{n=0}^{\infty} \left(\int_t^\infty v\int_0^\infty\Omega_k(vr)dF_n(r)dv\right)C_{n}^{(d-1)/2}(x).
\end{eqnarray*}

Using Fubini Theorem, we have 
\begin{eqnarray*}
	\int_t^\infty v\varphi(x,v)dv 
%	&=& \sum_{n=0}^{\infty}\left[\int_0^\infty \left(\int_t^\infty v\Omega_k(vr)dv\right)dF_n(r)\right] C_{n}^{(d-1)/2}(x)\\
	&=& 
	\sum_{n=0}^{\infty}\left[\int_0^\infty \left(\int_{tr}^\infty \frac{w}{r^2}\Omega_k(w)dw\right)dF_n(r)\right] C_{n}^{(d-1)/2}(x).
\end{eqnarray*}	

By \eqref{eq-int-Omega}, for $(x,t)\in[-1,1]\times[0,\infty)$, 
$$
\int_{tr}^\infty \frac{w}{r^2}\Omega_k(w)dw = \frac{(k-2)}{r^2}\Omega_{k-2}(tr).
$$

Hence, for $(x,t)\in[-1,1]\times[0,\infty)$,
\begin{eqnarray}\label{eq-int-numerador}
\int_t^\infty v\varphi(x,v)dv 
%&=&  
%(k-2)
%\sum_{n=0}^{\infty}\left[\int_0^\infty \Omega_{k-2}(tr)\frac{1}{r^2}dF_n(r)\right] C_{n}^{(d-1)/2}(x) \nonumber\\
&=&(k-2)
\sum_{n=0}^{\infty}g_n^1(t) C_{n}^{(d-1)/2}(x),
\end{eqnarray}	
where $g_n^{1}$ is defined in \eqref{eq-g_n(0-k)}.
%Now we have to evaluate \ana{(by hypothesis, it exists)} the following integral \ana{(in order to guarantee that it is positive)}
%$$
%\int_0^\infty v\varphi(x,v)dv.
%$$
%We have
In particular, 
\begin{eqnarray}\label{eq-int-dnominador}
\int_0^\infty v\varphi(1,v)dv  = 
%(k-2)\sum_{n=0}^{\infty}\left[\int_0^\infty \frac{1}{r^2}dF_n(r)\right] = 
(k-2)\sum_{n=0}^{\infty}g_n^{1}(0){C_{n}^{(d-1)/2}(1)},
\end{eqnarray}
%$$
%g_n^{1}(t)=\int_0^\infty \Omega_{k-2}(tr)\frac{1}{r^2}dF_n(r).
%$$
which is nonzero and finite. 
By \eqref{eq-int-numerador} and \eqref{eq-int-dnominador}, $I_2\varphi$ has the representation given in \eqref{eq-I2_0-k}.
\\

{\bf Induction step:} let assume the expression in \eqref{eq-I2_0-k} of $I_2^\kappa\varphi$ holds up to $\kappa$, and let us prove it holds for $I_2^{\kappa+1}\varphi$.
	
	We have
	\begin{eqnarray*}
		I_2^{\kappa+1}\varphi(x,t) &=& I_2(I_2^{\kappa}\varphi)(x,t) = 
		\frac1{\int_0^\infty vI_2^\kappa\varphi(1,v)dv} \int_t^\infty vI_2^\kappa\varphi(x,v)dv.
	\end{eqnarray*}
	
Note that the Hypothesis \eqref{h-int-r2-I2_0-k} guarantees that  $g_n^\kappa\in \mathcal P(\R^{k-2\kappa})$ and consequently the series in \eqref{eq-I2_0-k} converges absolutely and uniformly.

Using the induction hypothesis, integrating term by term, using Fubini theorem and \eqref{eq-int-Omega}, for $(x,t)\in[-1,1]\times[0,\infty)$, we obtain:
	\begin{eqnarray*}
		\int_t^\infty vI_2^\kappa\varphi(x,v)dv 
		&=&
	%	\frac1{\sum_{n=0}^{\infty}g_n^\kappa(0){C_{n}^{(d-1)/2}(1)}}
	\sum_{n=0}^{\infty} \left[\int_0^\infty \left(\int_t^\infty v\Omega_{k-2\kappa}(vr)dv\right)\frac{1}{r^{2\kappa}}dF_n(r) \right]  C_{n}^{(d-1)/2}(x)
		\\
&=&
(k-2\kappa-2)
		\sum_{n=0}^{\infty} \left[\int_0^\infty  \Omega_{k-2(\kappa+1)}(tr)\frac{1}{r^{2(\kappa+1)}}dF_n(r)\right]  C_{n}^{(d-1)/2}(x)	
		\\
		&=&
%		\frac{k-2\kappa-2}{\sum_{n=0}^{\infty}g_n^\kappa(0){C_{n}^{(d-1)/2}(1)}}
(k-2\kappa-2)
		\sum_{n=0}^{\infty} g_n^{\kappa+1}(t)  C_{n}^{(d-1)/2}(x).	
	\end{eqnarray*}

In particular,
	\begin{eqnarray*}
		\int_0^\infty vI_2^\kappa\varphi(1,v)dv =
(k-2\kappa-2)
\sum_{n=0}^{\infty} g_n^{\kappa+1}(0){C_{n}^{(d-1)/2}(1)},
	\end{eqnarray*}
which is nonzero and finite by \eqref{h-serie-I2_0-k}.	
	Therefore,	
	
	\begin{eqnarray*}
		I_2^{\kappa+1}\varphi(x,t) = 	\frac1{\sum_{n=0}^{\infty}g_n^{\kappa+1}(0){C_{n}^{(d-1)/2}(1)}}
		\sum_{n=0}^{\infty} g_n^{\kappa+1}(t)  C_{n}^{(d-1)/2}(x)
	\end{eqnarray*}	
%	with 
%	\begin{equation*}
%	g_n^{\kappa+1}(t) := \int_0^\infty \Omega_{k-2\kappa-2}(tr)\frac{1}{r^{2\kappa+2}}dF_n(r),
%	\end{equation*}
and \eqref{eq-I2_0-k} is proved.
	
Finally, given $\kappa\in\Z^*$, by \eqref{h-int-r2-I2_0-k} the $d$-Schoenberg functions $g_n^\kappa$ of     $I_2^\kappa\varphi$ belong to the class $\mathcal P(\R^{k-2\kappa})$ and together with  \eqref{h-serie-I2_0-k}  we can conclude  $0<\sum_{n=0}^{\infty} g_n^\kappa(0){C_{n}^{(d-1)/2}(1)}<\infty$. Therefore, Theorem 3.3 of \cite{porcu-berg} (see \eqref{eq-pd_SdxRk}) allows us to infer that  $I_2^\kappa\varphi$ belongs to $\mathcal P(\S^d\times\R^{k-2\kappa})$.

\end{proof}

\begin{proof}[{\bf Proof of Theorem \ref{t-I3k+H}}] 	We will prove \eqref{eq-I3k} by mathematical induction on $\kappa$.
	\\
	
%	{\bf Step $\kappa=0$:} Follows from \cite[Theorem 3.3]{porcu-berg} (see \eqref{eq-pd_SdxRk}).
%	\\

	{\bf Step $\kappa=1$:}
For each $(x,t)\in[-1,1]\times[0,\infty)$,
\begin{eqnarray*}
I_3^1\varphi(x,t) &=& %I_3\varphi(x,t)=
% \frac1{\int_0^\infty\int_{-1}^1 v\varphi(u,v)dudv} 
\int_t^\infty \int_{-1}^x v\varphi(u,v)dudv.
\end{eqnarray*}

Using \eqref{eq-pd_SdxRk} and  \eqref{eq-pd-Rk},
\begin{eqnarray*}
\int_t^\infty \int_{-1}^x v\varphi(u,v)dudv
% &=&  \frac1{\int_0^\infty\int_{-1}^1 v\varphi(u,v)dudv} 
% \int_t^\infty \int_{-1}^x v\sum_{n=0}^{\infty}f_{n}^d(v)C_{n}^{(d-1)/2}(u) dudv
%\\
&=& % \frac1{\int_0^\infty\int_{-1}^1 v\varphi(u,v)dudv} 
\int_t^\infty \int_{-1}^x v\sum_{n=0}^{\infty}\left(\int_0^\infty\Omega_k(vr)dF_n(r)\right)C_{n}^{(d-1)/2}(u) dudv.
\end{eqnarray*}

Integrating term by term and by  Fubbini Theorem, we have
\begin{eqnarray*}
	\int_t^\infty \int_{-1}^x v\varphi(u,v)dudv 	
%	&=& 
%	\sum_{n=0}^{\infty} \left[\int_t^\infty \int_{-1}^x v\left(\int_0^\infty\Omega_k(vr)dF_n(r)\right)C_{n}^{(d-1)/2}(u)du\right]\\
	&=& 
	\sum_{n=0}^{\infty} \left[\int_0^\infty \left(\int_t^\infty v\Omega_k(vr)dv\right)dF_n(r)\right] \left[\int_{-1}^xC_{n}^{(d-1)/2}(u)du\right].
\end{eqnarray*}

By \eqref{eq-int-Omega} and \eqref{eq-integ-Gegenb}, we obtain
\begin{eqnarray*}\label{eq-int_t_inf_-1_x}
\int_t^\infty \int_{-1}^x v\varphi(u,v)dudv	
&=& 
\sum_{n=0}^{\infty} \left[\int_0^\infty \frac{(k-2)}{r^2}\Omega_{k-2}(tr)dF_n(r)\right] 
\times
\\
&\times& \left[\frac1{\delta_{(d-1)/2}}(C_{n+1}^{(d-3)/2}(x)-C_{n+1}^{(d-3)/2}(-1))\right]
%\\
%&=& 
%\frac{(k-2)}{\delta_{(d-1)/2}}\sum_{n=0}^{\infty} g_n^1(t)  (C_{n+1}^{(d-3)/2}(x)-C_{n+1}^{(d-3)/2}(-1))
\end{eqnarray*}

Since
	$C_{n+1}^{(d-3)/2}(-1)=(-1)^{n+1}C_{n+1}^{(d-3)/2}(1)$, 
\begin{eqnarray*}
	\int_t^\infty \int_{-1}^x v\varphi(u,v)dudv	
&=& 
\frac{(k-2)}{\delta_{(d-1)/2}}
\left[
\left(\sum_{n=0}^{\infty} (-1)^nC_{n+1}^{(d-3)/2}(1)g_n^1(t)\right){C_0^{(d-3)/2}(x)} 
\right.
\\
&+&
\left.
\sum_{n=0}^{\infty}g_n^1(t)C_{n+1}^{(d-3)/2}(x)\right],
\end{eqnarray*}
where $g_n^1$ is given in \eqref{eq-g_n(0-k)}. 

Therefore,
\begin{eqnarray*}\label{eq-series_int_t_inf_-1_x}
I_3^1\varphi(x,t) &=& 
%\frac1{\gamma^{(d-3)/2,1}}
\sum_{n=0}^{\infty} h_n^1(t) C_{n}^{(d-3)/2}(x),
\end{eqnarray*}
where
\begin{eqnarray*}\label{eq-def-h_n}
h_n^1(t)= \left\{ 
\begin{array}{ll} 
\displaystyle	
%\frac{1}{2\sum_{j=0}^{\infty} g_{2j}^1(0)C_{2j+1}^{(d-3)/2}(1)} 
%\gamma^{d,1}
%\sum_{i=0}^{\infty} (-1)^i C_{i+1}^{(d-3)/2}(1) g_i^1(t), & n=0\\ 
%\\
%\gamma^{d,1}
\gamma^{d,k,1}
\sum_{i=0}^{\infty} (-1)^i \chi_i^{0,d,1} g_i^1(t), 
& n=0\\ \\
\displaystyle	
%\frac{1}{2\sum_{j=0}^{\infty} g_{2j}^1(0)C_{2j+1}^{(d-3)/2}(1)} 
%\gamma^{d,1}
\gamma^{d,k,1} g_{n-1}^1(t), & n\geq1,
\end{array}	
\right.
\end{eqnarray*}
where $\gamma^{d,k,1}=\frac{(k-2)}{\delta_{(d-1)/2}}>0$. Moreover
%$\gamma^{d,1}:=\left(2\sum_{j=0}^{\infty} g_{2j}^1(0)C_{2j+1}^{(d-3)/2}(1)\right)^{-1}$.
$\sum_{n=0}^{\infty} h_{n}^1(0)C_{n}^{(d-3)/2}(1)<\infty$ because, by \eqref{eq-chi} and \eqref{h-serie-gn-finit-I3k}-\eqref{h-serien-finit-I3k}, we have
\begin{eqnarray*}
\sum_{n=0}^{\infty} |h_{n}^1(0)C_{n}^{(d-3)/2}(1)|
&\leq&
%\sum_{i=0}^{\infty} |\chi_i^{0,d,1}| \, | g_i^1(0)| C_{0}^{(d-3)/2}(1)
%+
%\sum_{n=1}^{\infty} g_{n-1}^1(0) C_{n}^{(d-3)/2}(1)
%\\
%&\leq&
\varUpsilon^{0,d,1} C_{0}^{(d-3)/2}(1) \sum_{i=0}^{\infty} g_i^1(0) 
+
\sum_{n=1}^{\infty} g_{n-1}^1(0) C_{n}^{(d-3)/2}(1)<\infty.
\end{eqnarray*}

	{\bf Induction step:}  let assume the expression in \eqref{eq-I3k} of $I_3^\kappa\varphi$ holds up to $\kappa$, and let us prove it holds for $I_3^{\kappa+1}\varphi$.
	
	We have
	$$
	I_3^{\kappa+1}\varphi(x,t) = I_3(I_3^\kappa\varphi)(x,t) = 
%	\frac{1}{\int_0^\infty\int_{-1}^1 vI_3^\kappa\varphi(u,v)dudv}
	\int_t^\infty \int_{-1}^x vI_3^\kappa\varphi(u,v)dudv.
	$$
	
	Using the induction hypothesis, integrating term by term, using Fubini theorem, Equations \eqref{eq-g_n(0-k)}, \eqref{eq-int-Omega}, \eqref{eq-integ-Gegenb}, and making algebraic manipulations similar to the previous ones, for $(x,t)\in[-1,1]\times[0,\infty)$, we obtain
	\begin{eqnarray*}
		\int_t^\infty \int_{-1}^x vI_3^\kappa\varphi(u,v)dudv 
		&=& 
\gamma^{d,k,\kappa}
\sum_{n=0}^{\kappa-1}
\sum_{i=0}^{\infty}(-1)^{i}\chi_i^{n,d,\kappa}
\int_t^\infty  v  g_{i}^{\kappa}(v) dv 
\int_{-1}^x C_{n}^{(d-2\kappa-1)/2}(u)du
		\\
		&&
		+
%\frac{1}{\gamma^{(d-2\kappa-1)/2,\kappa}} 		
%\prod_{j=1}^\kappa\frac{k-2j}{\delta_{(d-2j+1)/2}}
\gamma^{d,k,\kappa}
\sum_{n=\kappa}^{\infty} 
\int_t^\infty v	g_{n-\kappa}^\kappa(v) dv
\int_{-1}^x  C_{n}^{(d-2\kappa-1)/2}(u)du
		\\
&=&
%\frac{1}{\gamma^{(d-2\kappa-1)/2,\kappa}}
\gamma^{d,k,\kappa}\frac{(k-2\kappa-2)}{\delta_{(d-2\kappa-1)/2}}	\times 
\\
&&
\hspace*{-10mm}\left[
\sum_{n=0}^{\kappa-1}
\sum_{i=0}^{\infty}(-1)^{i}\chi_i^{n,d,\kappa}
g_i^{\kappa+1}(t)
\left(
C_{n+1}^{(d-2\kappa-3)/2}(x)-C_{n+1}^{(d-2\kappa-3)/2}(-1)
\right)
\right. 
\\
&&
+ \left.
\sum_{n=\kappa}^{\infty}  
g_{n-\kappa}^{\kappa+1}(t)
\left(
C_{n+1}^{(d-2\kappa-3)/2}(x)-C_{n+1}^{(d-2\kappa-3)/2}(-1)
\right)
\right]
	\end{eqnarray*}
	
	Thus, as in the proof of Theorem \ref{t-I1},
\begin{eqnarray*}
	I_3^{\kappa+1}\varphi(x,t)  
	%		\int_t^\infty \int_{-1}^x vI_3^\kappa\varphi(u,v)dudv 
		&=&
%\frac{1}{\gamma^{(d-2\kappa-1)/2,\kappa}}\frac{(k-2\kappa-2)}{\delta_{(d-2\kappa-1)/2}}
\gamma^{d,k,\kappa+1}
\sum_{n=0}^{\infty} h_n^{\kappa+1}(t)C_{n}^{(d-2(\kappa+1)-1)/2}(x),
\end{eqnarray*}
where
\begin{eqnarray*}
h_n^{\kappa+1}(t):= \left\{ 
\begin{array}{lll} 
\displaystyle
\gamma^{d,k,\kappa+1}	
 \sum_{i=0}^{\infty}(-1)^{i}\chi_i^{n,d,\kappa+1}g_{i}^{\kappa+1}(t), & n=0,1,\ldots,\kappa\\ \\
 \gamma^{d,k,\kappa+1}
 g_{n-(\kappa+1)}^{\kappa+1}(t), & n\geq\kappa+1 \\ \\
\end{array}	
\right.
\end{eqnarray*}
with 
$\gamma^{d, k,\kappa+1}>0$, By \eqref{eq-chi} and \eqref{h-serie-gn-finit-I3k}-\eqref{h-serien-finit-I3k}, 
  $\sum_{n=0}^{\infty} h_{n}^{\kappa+1}(0)C_{n}^{(d-2(\kappa+1)-1)/2}(1)<\infty$
\end{proof}

\begin{proof}[{\bf Proof of Corollary \ref{c-tI3k+H}}]

We can proceed as in the proof of Corollary \ref{c-I1}  and rewrite $h_{n}^{\kappa}$, $n=0,1,\ldots,\kappa-1$,  as 
$$
h_{n}^{\kappa}(t) = \tilde h_{1,n}^\kappa(t) -\tilde h_{2,n}^\kappa(t),
$$
where 
\begin{eqnarray}\label{eq-til-h1k-h2k}
	&&
\tilde	h_{1,n}^\kappa(t) := %\frac{1}{\gamma^{(d-2\kappa-1)/2,\kappa}}
\gamma^{d, k,\kappa}
 \sum_{i=0}^{\infty}\chi_{2i}^{n,d,\kappa}g_{2i}^\kappa(t),\qquad \text{and} 
\nonumber	\\
	&&
\tilde	h_{2,n}^\kappa(t) := %\frac{1}{\gamma^{(d-2\kappa-1)/2,\kappa}}
\gamma^{d, k,\kappa}
\sum_{i=0}^{\infty}\chi_{2i+1}^{n,d,\kappa}g_{2i+1}^\kappa(t).
\end{eqnarray}
Define the bounded function $H^\kappa$ on $[-1,1]\times[0,\infty)$ by 
$$
H^\kappa(x,t):=  \sum_{n=0}^{\kappa-1}
\tilde h_{2,n}^{\kappa}(t)C_{n}^{(d-2\kappa-1)/2}(x) -  
\tilde h_{1,0}^{\kappa}(t)C_{0}^{(d-2\kappa-1)/2}(x) 
$$

By Remark \ref{r-chi-positivo}, it is clear that $\tilde h_{1,n}^\kappa\in\mathcal P(\R^{k})$, $n=1,2,\ldots,\kappa-1$, and also $h_{n}^{\kappa}\in\mathcal P(\R^{k})$ for $n\geq\kappa$. Therefore,
$$
H^\kappa(x,t)+I_3^\kappa\varphi(x,t) =  
\sum_{n=1}^{\kappa-1}\tilde h_{1,n}^\kappa(t) C_{n}^{(d-2\kappa-1)/2}(x) 
+
%\frac1{\gamma^{(d-2\kappa-1)/2,\kappa}}
\sum_{n=\kappa}^{\infty} h_{n}^{\kappa}(t)C_{n}^{(d-2\kappa-1)/2}(x)
$$
has an expansion as \eqref{eq-pd_SdxRk} with the series uniformly convergent on $[-1,1]\times[0,\infty)$
due to Theorem \ref{t-I3k+H}.
By Theorem 3.3 of \cite{porcu-berg} (see \eqref{eq-pd_SdxRk}), we can conclude that the function $H^\kappa+I_1^\kappa\varphi$ belongs to the class $\mathcal P(\S^{d-2\kappa}\times\R^{k})$.
\end{proof}

\begin{proof}[{\bf Proof of Theorem \ref{t-I3k-pd-A-B}    }] 
	
As in the proof of Theorem \ref{t-I1k-AB},  for any  constants $A$ and $B$, by \eqref{eq-g_n(0-k)},
	\begin{equation}\label{eq-h1-h2-I3}
		A\tilde h_{1,n}^\kappa(t)  - B\tilde h_{2,n}^\kappa(t) 
		=
%\frac{1}{\gamma^{(d-2\kappa-1)/2,\kappa}}
		\sum_{i=0}^{\infty} \int_0^\infty  \Omega_{k-2\kappa}(tr)
		\left(A\chi_{2i}^{n,d,\kappa}\frac1{r^{2\kappa}}dF_{2i}(r) - B\chi_{2i+1}^{n,d,\kappa}\frac1{r^{2\kappa}}dF_{2i+1}(r)\right).
	\end{equation}	
	
	Since, by \eqref{eq-chi},
	$$
	\int_0^\infty A|\chi_{i}^{n,d,\kappa}|\frac1{r^{2\kappa}}dF_{i}(r) 
	\leq
	A\varUpsilon^{n,d,\kappa}C_i^{(d-1)/2}(1)\int_0^\infty \frac1{r^{2\kappa}} dF_{i}(r),
	$$
	we have
	$$
	\int_0^\infty\left(A\chi_{2i}^{n,d,\kappa}\frac1{r^{2\kappa}}dF_{2i}(r) - B\chi_{2i+1}^{n,d,\kappa}\frac1{r^{2\kappa}}dF_{2i+1}(r)\right)<\infty
	$$
	and by Hypothesis \eqref{h-serie-gn-Cn-finit-I3k} the series  in \eqref{eq-h1-h2-I3} converges absolutely and uniformly on $[0,\infty)$. 
	
	Thus,
	\begin{eqnarray*}
		A\tilde h_{1,n}^\kappa(t)  - B\tilde h_{2,n}^\kappa(t) 
		&=&
%\frac{1}{\gamma^{(d-2\kappa-1)/2,\kappa}}
		\int_0^\infty  \Omega_k(tr)
		d\left(	
		\sum_{i=0}^{\infty} 
		A\chi_{2i}^{n,d,\kappa}\frac1{r^{2\kappa}}F_{2i}(r) - B\chi_{2i+1}^{n,d,\kappa}\frac1{r^{2\kappa}}F_{2i+1}(r)\right).
	\end{eqnarray*}	
	
	By \eqref{eq-cond-medida_F_n},  the series $\sum_{i=0}^{\infty} 
	\chi_{2i}^{n,d,\kappa}\frac1{r^{2\kappa}}F_{2i}$ and $\sum_{i=0}^{\infty}  \chi_{2i+1}^{n,d,\kappa}\frac1{r^{2\kappa}}F_{2i+1}$ are uniformly bounded on $[0,\infty)$. 
	Then we can choose $A^n,B^n$ such that
	$\sum_{i=0}^{\infty} 
	A^n\chi_{2i}^{n,d,\kappa}\frac1{r^{2\kappa}}F_{2i} - B^n\chi_{2i+1}^{n,d,\kappa}\frac1{r^{2\kappa}}F_{2i+1}$ is non negative which allows us conclude that  $A^n\tilde h_{1,n}^\kappa-B^n\tilde h_{2,n}^\kappa\in\mathcal P(\R^{k})$.
The uniform convergence of the series \eqref{eq-I3kAB} follows by Theorem \ref{t-I3k+H} and the result  by Theorem 3.3 of \cite{porcu-berg} (see \eqref{eq-pd_SdxRk}).	
	
\end{proof}

\bibliographystyle{apalike}
%\bibliographystyle{mywiley}
%\bibliography{bib-new}{}

%\bibliography{/../AnaPaula/BIB_ana_latex/bibAna}

\bibliography{bibwalksproduct}

\end{document}